\documentclass[10pt,a4paper]{article}

\usepackage{amsfonts}
\usepackage{amsmath}
\usepackage{amssymb}
\usepackage{graphicx}
\usepackage{harpoon}
\usepackage[
    colorlinks,
    citecolor=darkgray,
    filecolor=darkgray,
    linkcolor=darkgray,
    urlcolor=darkgray]{hyperref}
\usepackage{tikz}
\usepackage{xcolor}

\author{R\'emy Sigrist}
\date{\today}
\title{The Kochawave curve, a variant of the Koch curve}

\newcommand{\quotes}[1]{``#1''}
\newcommand{\twodots}{\mathinner {\ldotp \ldotp}}

\begin{document}

\maketitle

\begin{abstract}

This paper presents the construction and some properties
of the Kochawave curve, a simple and asymmetrical variant of the Koch curve.

\end{abstract}

\begin{figure}[b]
    \begin{center}
    \includegraphics[height=1.5cm, keepaspectratio]{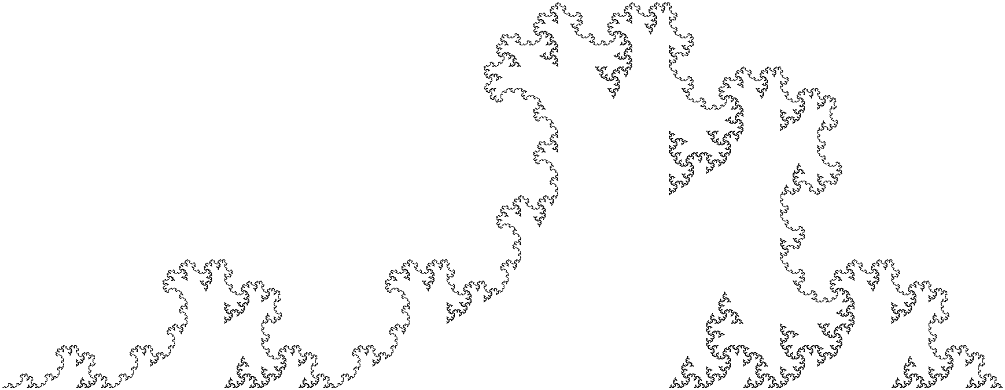}
    \end{center}
\end{figure}

\newpage

\tableofcontents

\newpage

\section{Introduction}

The Koch curve~\cite{wiki} is a well-known fractal described in 1904
by the Swedish mathematician Helge von Koch~\cite{koch}.

Since then, many variants of the Koch curve have been described
in the literature, such as the Ces\`aro fractal~\cite{weisstein}
and the quadratic Koch curve~\cite{ferreol},
just to name two (see Li~\cite{li}, Rani~\cite{rani}, Ventrella~\cite{ventrella} and Wikipedia~\cite{wiki}
for other variants).

We present here another variant of the Koch curve:
the Kochawave curve, mentioned and named by Wahl~\cite{wahl}.
The Kochawave curve has a simple and asymmetrical structure,
and as its name suggests, looks like a wave.

\section{Construction}

The Kochawave curve is a fractal curve.
We describe four ways to construct it.

\subsection{Construction with segments}
\label{cons-segments}

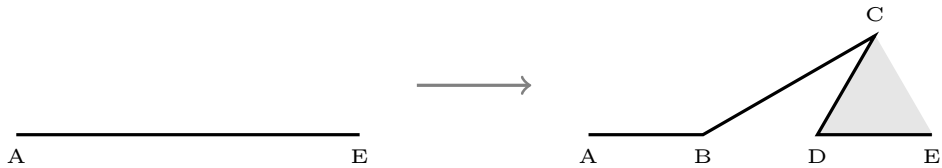
\begin{figure}[!ht]
    \resizebox{13cm}{!}
    {
    \begin{tikzpicture}
    	\begin{scope}[yscale=.87,xslant=.5]
    		\draw[thick] (0,0) node[below] {\tiny A} -- (3,0) node[below] {\tiny E};
    		\draw[->,thick,gray] (3.25,0.5) -> (4.25,0.5);
    		\draw[fill,gray!20] (7,0) -- (8,0) -- (7,1);
    		\draw[thick] (5,0) node[below] {\tiny A} -- (6,0) node[below] {\tiny B} -- (7,1) node[above] {\tiny C} -- (7,0) node[below] {\tiny D} -- (8,0) node[below] {\tiny E};
    	\end{scope}
    \end{tikzpicture}
    }
    \caption{The construction rule for the Kochawave curve}
    \label{fig:seg-one-step}
\end{figure}

The Kochawave curve can be constructed starting from one unit segment,
and iteratively replacing each line segment $\overrightarrow{AE}$
with four line segments $\overrightarrow{AB}$, $\overrightarrow{BC}$, $\overrightarrow{CD}$, $\overrightarrow{DE}$ such that
the points $A$, $B$, $D$ and $E$ are colinear and equally spaced (and appear in that order) and $DEC$ forms an equilateral triangle
(see figure~\ref{fig:seg-one-step}).

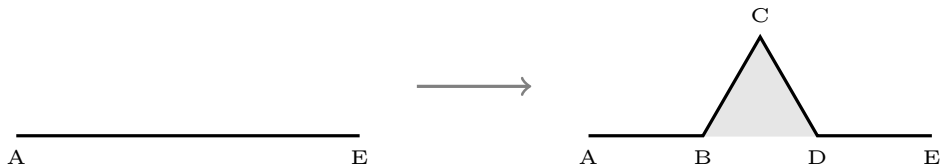
\begin{figure}[!ht]
    \resizebox{13cm}{!}
    {
    \begin{tikzpicture}
    	\begin{scope}[yscale=.87,xslant=.5]
    		\draw[thick] (0,0) node[below] {\tiny A} -- (3,0) node[below] {\tiny E};
    		\draw[->,thick,gray] (3.25,0.5) -> (4.25,0.5);
    		\draw[fill,gray!20] (6,0) -- (7,0) -- (6,1);
    		\draw[thick] (5,0) node[below] {\tiny A} -- (6,0) node[below] {\tiny B} -- (6,1) node[above] {\tiny C} -- (7,0) node[below] {\tiny D} -- (8,0) node[below] {\tiny E};
    	\end{scope}
    \end{tikzpicture}
    }
    \caption{The construction rule for the Koch curve}
    \label{fig:org-one-step}
\end{figure}

The Koch curve has a similar construction rule, except that there $BDC$ forms an equilateral triangle (see figure~\ref{fig:org-one-step}).

\begin{figure}[!ht]
    \resizebox{13cm}{!}
    {
    \begin{tikzpicture}
    	\begin{scope}[yscale=.87,xslant=.5]
    		\draw[line width=5pt] (0,0) -- (27,0) -- (54,27) -- (54,0) -- (81,0);
    		\draw[line width=5pt] (91,0) -- (100,0) -- (109,9) -- (109,0) -- (118,0) -- (127,9) -- (127,36) -- (136,18) -- (145,27) -- (145,18) -- (154,0) -- (145,9) -- (145,0) -- (154,0) -- (163,9) -- (163,0) -- (172,0);
    		\draw[line width=5pt] (182,0) -- (185,0) -- (188,3) -- (188,0) -- (191,0) -- (194,3) -- (194,12) -- (197,6) -- (200,9) -- (200,6) -- (203,0) -- (200,3) -- (200,0) -- (203,0) -- (206,3) -- (206,0) -- (209,0) -- (212,3) -- (212,12) -- (215,6) -- (218,9) -- (218,18) -- (209,36) -- (218,27) -- (218,36) -- (221,30) -- (230,21) -- (224,24) -- (227,18) -- (230,21) -- (230,30) -- (233,24) -- (236,27) -- (236,24) -- (239,18) -- (236,21) -- (236,18) -- (239,12) -- (248,3) -- (242,6) -- (245,0) -- (242,3) -- (236,6) -- (239,6) -- (236,9) -- (236,6) -- (239,0) -- (236,3) -- (236,0) -- (239,0) -- (242,3) -- (242,0) -- (245,0) -- (248,3) -- (248,12) -- (251,6) -- (254,9) -- (254,6) -- (257,0) -- (254,3) -- (254,0) -- (257,0) -- (260,3) -- (260,0) -- (263,0);
    		\draw[line width=5pt] (273,0) -- (274,0) -- (275,1) -- (275,0) -- (276,0) -- (277,1) -- (277,4) -- (278,2) -- (279,3) -- (279,2) -- (280,0) -- (279,1) -- (279,0) -- (280,0) -- (281,1) -- (281,0) -- (282,0) -- (283,1) -- (283,4) -- (284,2) -- (285,3) -- (285,6) -- (282,12) -- (285,9) -- (285,12) -- (286,10) -- (289,7) -- (287,8) -- (288,6) -- (289,7) -- (289,10) -- (290,8) -- (291,9) -- (291,8) -- (292,6) -- (291,7) -- (291,6) -- (292,4) -- (295,1) -- (293,2) -- (294,0) -- (293,1) -- (291,2) -- (292,2) -- (291,3) -- (291,2) -- (292,0) -- (291,1) -- (291,0) -- (292,0) -- (293,1) -- (293,0) -- (294,0) -- (295,1) -- (295,4) -- (296,2) -- (297,3) -- (297,2) -- (298,0) -- (297,1) -- (297,0) -- (298,0) -- (299,1) -- (299,0) -- (300,0) -- (301,1) -- (301,4) -- (302,2) -- (303,3) -- (303,6) -- (300,12) -- (303,9) -- (303,12) -- (304,10) -- (307,7) -- (305,8) -- (306,6) -- (307,7) -- (307,10) -- (308,8) -- (309,9) -- (309,12) -- (306,18) -- (309,15) -- (309,18) -- (306,24) -- (297,33) -- (303,30) -- (300,36) -- (303,33) -- (309,30) -- (306,30) -- (309,27) -- (309,30) -- (306,36) -- (309,33) -- (309,36) -- (310,34) -- (313,31) -- (311,32) -- (312,30) -- (315,27) -- (321,24) -- (318,24) -- (321,21) -- (319,22) -- (316,22) -- (317,23) -- (315,24) -- (316,22) -- (319,19) -- (317,20) -- (318,18) -- (319,19) -- (319,22) -- (320,20) -- (321,21) -- (321,24) -- (318,30) -- (321,27) -- (321,30) -- (322,28) -- (325,25) -- (323,26) -- (324,24) -- (325,25) -- (325,28) -- (326,26) -- (327,27) -- (327,26) -- (328,24) -- (327,25) -- (327,24) -- (328,22) -- (331,19) -- (329,20) -- (330,18) -- (329,19) -- (327,20) -- (328,20) -- (327,21) -- (327,20) -- (328,18) -- (327,19) -- (327,18) -- (328,16) -- (331,13) -- (329,14) -- (330,12) -- (333,9) -- (339,6) -- (336,6) -- (339,3) -- (337,4) -- (334,4) -- (335,5) -- (333,6) -- (334,4) -- (337,1) -- (335,2) -- (336,0) -- (335,1) -- (333,2) -- (334,2) -- (333,3) -- (331,4) -- (328,4) -- (329,5) -- (327,6) -- (328,6) -- (329,7) -- (329,6) -- (330,6) -- (329,7) -- (327,8) -- (328,8) -- (327,9) -- (327,8) -- (328,6) -- (327,7) -- (327,6) -- (328,4) -- (331,1) -- (329,2) -- (330,0) -- (329,1) -- (327,2) -- (328,2) -- (327,3) -- (327,2) -- (328,0) -- (327,1) -- (327,0) -- (328,0) -- (329,1) -- (329,0) -- (330,0) -- (331,1) -- (331,4) -- (332,2) -- (333,3) -- (333,2) -- (334,0) -- (333,1) -- (333,0) -- (334,0) -- (335,1) -- (335,0) -- (336,0) -- (337,1) -- (337,4) -- (338,2) -- (339,3) -- (339,6) -- (336,12) -- (339,9) -- (339,12) -- (340,10) -- (343,7) -- (341,8) -- (342,6) -- (343,7) -- (343,10) -- (344,8) -- (345,9) -- (345,8) -- (346,6) -- (345,7) -- (345,6) -- (346,4) -- (349,1) -- (347,2) -- (348,0) -- (347,1) -- (345,2) -- (346,2) -- (345,3) -- (345,2) -- (346,0) -- (345,1) -- (345,0) -- (346,0) -- (347,1) -- (347,0) -- (348,0) -- (349,1) -- (349,4) -- (350,2) -- (351,3) -- (351,2) -- (352,0) -- (351,1) -- (351,0) -- (352,0) -- (353,1) -- (353,0) -- (354,0);
    	\end{scope}
    \end{tikzpicture}
    }
    \caption{The first four iterations of the Kochawave curve}
    \label{fig:seg-steps}
\end{figure}
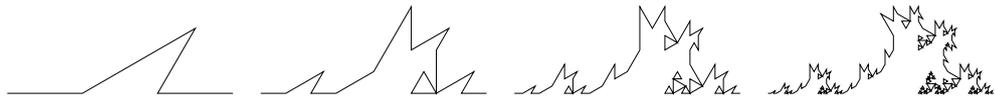

The first four iterations of the construction of the Kochawave curve are depicted in figure~\ref{fig:seg-steps}.

\begin{figure}[!ht]
  \includegraphics[width=\linewidth]{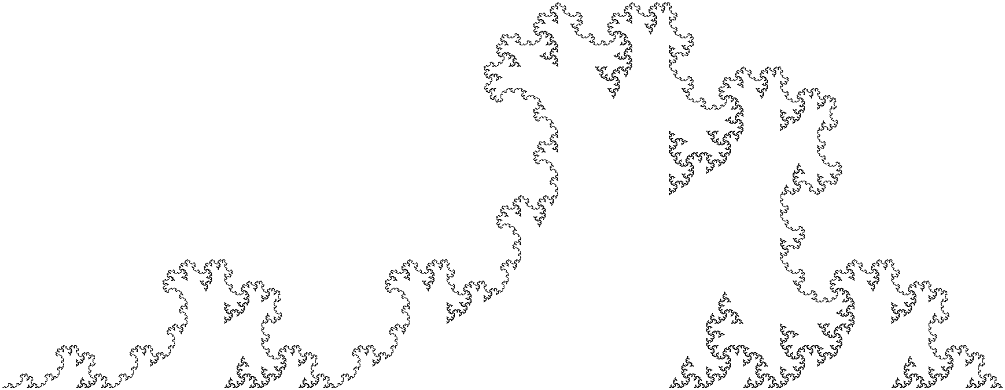}
  \caption{The Kochawave curve}
  \label{fig:kochawave}
\end{figure}

The Kochawave curve is the limiting figure as the iterations continue indefinitely (see figure~\ref{fig:kochawave}).

\clearpage

\subsection{Construction with triangles}
\label{cons-triangles}

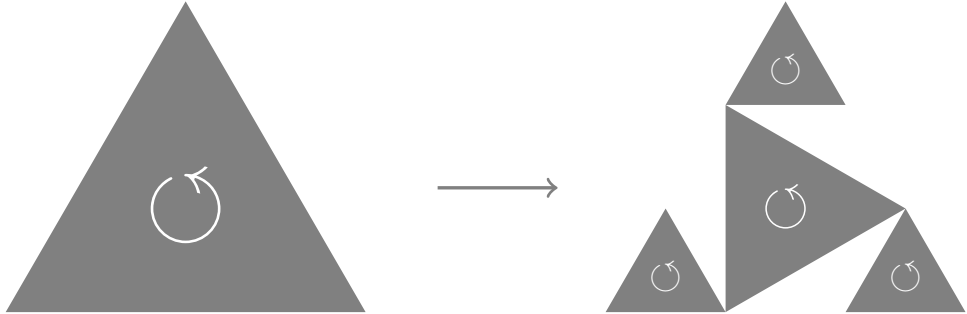
\begin{figure}[!ht]
    \resizebox{13cm}{!}
    {
    \begin{tikzpicture}
    	\begin{scope}[yscale=.87,xslant=.5]
    		\coordinate (A) at (0,3);
    		\coordinate (B) at (0,0);
    		\coordinate (C) at (3,0);
    		\fill[gray] (A) -- (B) -- (C);

    		\node[white] at (barycentric cs:A=1,B=1,C=1) {\Huge $\circlearrowleft$};

    		\draw[->,thick,gray] (3,1.2) -> (4,1.2);

    		\coordinate (a) at (5,3);
    		\coordinate (b) at (5,2);
    		\coordinate (c) at (6,2);
    		\coordinate (d) at (5,1);
    		\coordinate (e) at (7,1);
    		\coordinate (f) at (5,0);
    		\coordinate (g) at (6,0);
    		\coordinate (h) at (7,0);
    		\coordinate (i) at (8,0);
    		\fill[gray] (a) -- (b) -- (g) -- (d) -- (f) -- (g) -- (e) -- (h) -- (i) -- (e) -- (b) -- (c);

    		\node[white] at (barycentric cs:a=1,b=1,c=1) {$\circlearrowleft$};
    		\node[white] at (barycentric cs:d=1,f=1,g=1) {$\circlearrowleft$};
    		\node[white] at (barycentric cs:e=1,h=1,i=1) {$\circlearrowleft$};
    		\node[white] at (barycentric cs:b=1,e=1,g=1) {\Large $\circlearrowleft$};
    	\end{scope}
    \end{tikzpicture}
    }
    \caption{The construction rule for the three entangled Kochawave curves}
    \label{fig:tri-one-step}
\end{figure}

We can construct three entangled Kochawave curves starting from one unit sided equilateral triangle,
and iteratively replacing each equilateral triangle, say of side length $w$,
with four equilateral triangles, of side lengths $w/3$, $w/3$, $w/3$ and $w/\sqrt{3}$, respectively,
with the same orientation as in the original triangle, as depicted in figure~\ref{fig:tri-one-step}.

\begin{figure}[!ht]
    \resizebox{13cm}{!}
    {
    \begin{tikzpicture}
    	\begin{scope}[yscale=.87,xslant=.5]
    		\fill[gray] (0,0) -- (27,0) -- (54,27) -- (54,0) -- (81,0) -- (54,27) -- (0,54) -- (27,54) -- (0,81) -- (0,54) -- (27,0) -- (0,27) -- (0,0);
    		\fill[gray] (91,0) -- (100,0) -- (109,9) -- (109,0) -- (118,0) -- (127,9) -- (127,36) -- (136,18) -- (145,27) -- (145,18) -- (154,0) -- (145,9) -- (145,0) -- (154,0) -- (163,9) -- (163,0) -- (172,0) -- (163,9) -- (145,18) -- (154,18) -- (145,27) -- (127,36) -- (100,36) -- (109,45) -- (91,54) -- (100,54) -- (109,63) -- (109,54) -- (118,54) -- (109,63) -- (91,72) -- (100,72) -- (91,81) -- (91,72) -- (100,54) -- (91,63) -- (91,54) -- (100,36) -- (127,9) -- (109,18) -- (118,0) -- (109,9) -- (91,18) -- (100,18) -- (91,27) -- (91,18) -- (100,0) -- (91,9) -- (91,0);
    		\fill[gray] (182,0) -- (185,0) -- (188,3) -- (188,0) -- (191,0) -- (194,3) -- (194,12) -- (197,6) -- (200,9) -- (200,6) -- (203,0) -- (200,3) -- (200,0) -- (203,0) -- (206,3) -- (206,0) -- (209,0) -- (212,3) -- (212,12) -- (215,6) -- (218,9) -- (218,18) -- (209,36) -- (218,27) -- (218,36) -- (221,30) -- (230,21) -- (224,24) -- (227,18) -- (230,21) -- (230,30) -- (233,24) -- (236,27) -- (236,24) -- (239,18) -- (236,21) -- (236,18) -- (239,12) -- (248,3) -- (242,6) -- (245,0) -- (242,3) -- (236,6) -- (239,6) -- (236,9) -- (236,6) -- (239,0) -- (236,3) -- (236,0) -- (239,0) -- (242,3) -- (242,0) -- (245,0) -- (248,3) -- (248,12) -- (251,6) -- (254,9) -- (254,6) -- (257,0) -- (254,3) -- (254,0) -- (257,0) -- (260,3) -- (260,0) -- (263,0) -- (260,3) -- (254,6) -- (257,6) -- (254,9) -- (248,12) -- (239,12) -- (242,15) -- (236,18) -- (239,18) -- (242,21) -- (242,18) -- (245,18) -- (242,21) -- (236,24) -- (239,24) -- (236,27) -- (230,30) -- (221,30) -- (224,33) -- (218,36) -- (209,36) -- (200,27) -- (200,36) -- (191,36) -- (194,39) -- (194,48) -- (197,42) -- (200,45) -- (194,48) -- (185,48) -- (188,51) -- (182,54) -- (185,54) -- (188,57) -- (188,54) -- (191,54) -- (194,57) -- (194,66) -- (197,60) -- (200,63) -- (200,60) -- (203,54) -- (200,57) -- (200,54) -- (203,54) -- (206,57) -- (206,54) -- (209,54) -- (206,57) -- (200,60) -- (203,60) -- (200,63) -- (194,66) -- (185,66) -- (188,69) -- (182,72) -- (185,72) -- (188,75) -- (188,72) -- (191,72) -- (188,75) -- (182,78) -- (185,78) -- (182,81) -- (182,78) -- (185,72) -- (182,75) -- (182,72) -- (185,66) -- (194,57) -- (188,60) -- (191,54) -- (188,57) -- (182,60) -- (185,60) -- (182,63) -- (182,60) -- (185,54) -- (182,57) -- (182,54) -- (185,48) -- (194,39) -- (188,42) -- (191,36) -- (200,27) -- (218,18) -- (209,18) -- (218,9) -- (212,12) -- (203,12) -- (206,15) -- (200,18) -- (203,12) -- (212,3) -- (206,6) -- (209,0) -- (206,3) -- (200,6) -- (203,6) -- (200,9) -- (194,12) -- (185,12) -- (188,15) -- (182,18) -- (185,18) -- (188,21) -- (188,18) -- (191,18) -- (188,21) -- (182,24) -- (185,24) -- (182,27) -- (182,24) -- (185,18) -- (182,21) -- (182,18) -- (185,12) -- (194,3) -- (188,6) -- (191,0) -- (188,3) -- (182,6) -- (185,6) -- (182,9) -- (182,6) -- (185,0) -- (182,3) -- (182,0);
    		\fill[gray] (273,0) -- (274,0) -- (275,1) -- (275,0) -- (276,0) -- (277,1) -- (277,4) -- (278,2) -- (279,3) -- (279,2) -- (280,0) -- (279,1) -- (279,0) -- (280,0) -- (281,1) -- (281,0) -- (282,0) -- (283,1) -- (283,4) -- (284,2) -- (285,3) -- (285,6) -- (282,12) -- (285,9) -- (285,12) -- (286,10) -- (289,7) -- (287,8) -- (288,6) -- (289,7) -- (289,10) -- (290,8) -- (291,9) -- (291,8) -- (292,6) -- (291,7) -- (291,6) -- (292,4) -- (295,1) -- (293,2) -- (294,0) -- (293,1) -- (291,2) -- (292,2) -- (291,3) -- (291,2) -- (292,0) -- (291,1) -- (291,0) -- (292,0) -- (293,1) -- (293,0) -- (294,0) -- (295,1) -- (295,4) -- (296,2) -- (297,3) -- (297,2) -- (298,0) -- (297,1) -- (297,0) -- (298,0) -- (299,1) -- (299,0) -- (300,0) -- (301,1) -- (301,4) -- (302,2) -- (303,3) -- (303,6) -- (300,12) -- (303,9) -- (303,12) -- (304,10) -- (307,7) -- (305,8) -- (306,6) -- (307,7) -- (307,10) -- (308,8) -- (309,9) -- (309,12) -- (306,18) -- (309,15) -- (309,18) -- (306,24) -- (297,33) -- (303,30) -- (300,36) -- (303,33) -- (309,30) -- (306,30) -- (309,27) -- (309,30) -- (306,36) -- (309,33) -- (309,36) -- (310,34) -- (313,31) -- (311,32) -- (312,30) -- (315,27) -- (321,24) -- (318,24) -- (321,21) -- (319,22) -- (316,22) -- (317,23) -- (315,24) -- (316,22) -- (319,19) -- (317,20) -- (318,18) -- (319,19) -- (319,22) -- (320,20) -- (321,21) -- (321,24) -- (318,30) -- (321,27) -- (321,30) -- (322,28) -- (325,25) -- (323,26) -- (324,24) -- (325,25) -- (325,28) -- (326,26) -- (327,27) -- (327,26) -- (328,24) -- (327,25) -- (327,24) -- (328,22) -- (331,19) -- (329,20) -- (330,18) -- (329,19) -- (327,20) -- (328,20) -- (327,21) -- (327,20) -- (328,18) -- (327,19) -- (327,18) -- (328,16) -- (331,13) -- (329,14) -- (330,12) -- (333,9) -- (339,6) -- (336,6) -- (339,3) -- (337,4) -- (334,4) -- (335,5) -- (333,6) -- (334,4) -- (337,1) -- (335,2) -- (336,0) -- (335,1) -- (333,2) -- (334,2) -- (333,3) -- (331,4) -- (328,4) -- (329,5) -- (327,6) -- (328,6) -- (329,7) -- (329,6) -- (330,6) -- (329,7) -- (327,8) -- (328,8) -- (327,9) -- (327,8) -- (328,6) -- (327,7) -- (327,6) -- (328,4) -- (331,1) -- (329,2) -- (330,0) -- (329,1) -- (327,2) -- (328,2) -- (327,3) -- (327,2) -- (328,0) -- (327,1) -- (327,0) -- (328,0) -- (329,1) -- (329,0) -- (330,0) -- (331,1) -- (331,4) -- (332,2) -- (333,3) -- (333,2) -- (334,0) -- (333,1) -- (333,0) -- (334,0) -- (335,1) -- (335,0) -- (336,0) -- (337,1) -- (337,4) -- (338,2) -- (339,3) -- (339,6) -- (336,12) -- (339,9) -- (339,12) -- (340,10) -- (343,7) -- (341,8) -- (342,6) -- (343,7) -- (343,10) -- (344,8) -- (345,9) -- (345,8) -- (346,6) -- (345,7) -- (345,6) -- (346,4) -- (349,1) -- (347,2) -- (348,0) -- (347,1) -- (345,2) -- (346,2) -- (345,3) -- (345,2) -- (346,0) -- (345,1) -- (345,0) -- (346,0) -- (347,1) -- (347,0) -- (348,0) -- (349,1) -- (349,4) -- (350,2) -- (351,3) -- (351,2) -- (352,0) -- (351,1) -- (351,0) -- (352,0) -- (353,1) -- (353,0) -- (354,0) -- (353,1) -- (351,2) -- (352,2) -- (351,3) -- (349,4) -- (346,4) -- (347,5) -- (345,6) -- (346,6) -- (347,7) -- (347,6) -- (348,6) -- (347,7) -- (345,8) -- (346,8) -- (345,9) -- (343,10) -- (340,10) -- (341,11) -- (339,12) -- (336,12) -- (333,9) -- (333,12) -- (330,12) -- (331,13) -- (331,16) -- (332,14) -- (333,15) -- (331,16) -- (328,16) -- (329,17) -- (327,18) -- (328,18) -- (329,19) -- (329,18) -- (330,18) -- (331,19) -- (331,22) -- (332,20) -- (333,21) -- (333,20) -- (334,18) -- (333,19) -- (333,18) -- (334,18) -- (335,19) -- (335,18) -- (336,18) -- (335,19) -- (333,20) -- (334,20) -- (333,21) -- (331,22) -- (328,22) -- (329,23) -- (327,24) -- (328,24) -- (329,25) -- (329,24) -- (330,24) -- (329,25) -- (327,26) -- (328,26) -- (327,27) -- (325,28) -- (322,28) -- (323,29) -- (321,30) -- (318,30) -- (315,27) -- (315,30) -- (312,30) -- (313,31) -- (313,34) -- (314,32) -- (315,33) -- (313,34) -- (310,34) -- (311,35) -- (309,36) -- (306,36) -- (303,33) -- (303,36) -- (300,36) -- (297,33) -- (297,24) -- (294,30) -- (291,27) -- (291,30) -- (288,36) -- (291,33) -- (291,36) -- (288,36) -- (285,33) -- (285,36) -- (282,36) -- (283,37) -- (283,40) -- (284,38) -- (285,39) -- (285,42) -- (282,48) -- (285,45) -- (285,48) -- (286,46) -- (289,43) -- (287,44) -- (288,42) -- (289,43) -- (289,46) -- (290,44) -- (291,45) -- (289,46) -- (286,46) -- (287,47) -- (285,48) -- (282,48) -- (279,45) -- (279,48) -- (276,48) -- (277,49) -- (277,52) -- (278,50) -- (279,51) -- (277,52) -- (274,52) -- (275,53) -- (273,54) -- (274,54) -- (275,55) -- (275,54) -- (276,54) -- (277,55) -- (277,58) -- (278,56) -- (279,57) -- (279,56) -- (280,54) -- (279,55) -- (279,54) -- (280,54) -- (281,55) -- (281,54) -- (282,54) -- (283,55) -- (283,58) -- (284,56) -- (285,57) -- (285,60) -- (282,66) -- (285,63) -- (285,66) -- (286,64) -- (289,61) -- (287,62) -- (288,60) -- (289,61) -- (289,64) -- (290,62) -- (291,63) -- (291,62) -- (292,60) -- (291,61) -- (291,60) -- (292,58) -- (295,55) -- (293,56) -- (294,54) -- (293,55) -- (291,56) -- (292,56) -- (291,57) -- (291,56) -- (292,54) -- (291,55) -- (291,54) -- (292,54) -- (293,55) -- (293,54) -- (294,54) -- (295,55) -- (295,58) -- (296,56) -- (297,57) -- (297,56) -- (298,54) -- (297,55) -- (297,54) -- (298,54) -- (299,55) -- (299,54) -- (300,54) -- (299,55) -- (297,56) -- (298,56) -- (297,57) -- (295,58) -- (292,58) -- (293,59) -- (291,60) -- (292,60) -- (293,61) -- (293,60) -- (294,60) -- (293,61) -- (291,62) -- (292,62) -- (291,63) -- (289,64) -- (286,64) -- (287,65) -- (285,66) -- (282,66) -- (279,63) -- (279,66) -- (276,66) -- (277,67) -- (277,70) -- (278,68) -- (279,69) -- (277,70) -- (274,70) -- (275,71) -- (273,72) -- (274,72) -- (275,73) -- (275,72) -- (276,72) -- (277,73) -- (277,76) -- (278,74) -- (279,75) -- (279,74) -- (280,72) -- (279,73) -- (279,72) -- (280,72) -- (281,73) -- (281,72) -- (282,72) -- (281,73) -- (279,74) -- (280,74) -- (279,75) -- (277,76) -- (274,76) -- (275,77) -- (273,78) -- (274,78) -- (275,79) -- (275,78) -- (276,78) -- (275,79) -- (273,80) -- (274,80) -- (273,81) -- (273,80) -- (274,78) -- (273,79) -- (273,78) -- (274,76) -- (277,73) -- (275,74) -- (276,72) -- (275,73) -- (273,74) -- (274,74) -- (273,75) -- (273,74) -- (274,72) -- (273,73) -- (273,72) -- (274,70) -- (277,67) -- (275,68) -- (276,66) -- (279,63) -- (285,60) -- (282,60) -- (285,57) -- (283,58) -- (280,58) -- (281,59) -- (279,60) -- (280,58) -- (283,55) -- (281,56) -- (282,54) -- (281,55) -- (279,56) -- (280,56) -- (279,57) -- (277,58) -- (274,58) -- (275,59) -- (273,60) -- (274,60) -- (275,61) -- (275,60) -- (276,60) -- (275,61) -- (273,62) -- (274,62) -- (273,63) -- (273,62) -- (274,60) -- (273,61) -- (273,60) -- (274,58) -- (277,55) -- (275,56) -- (276,54) -- (275,55) -- (273,56) -- (274,56) -- (273,57) -- (273,56) -- (274,54) -- (273,55) -- (273,54) -- (274,52) -- (277,49) -- (275,50) -- (276,48) -- (279,45) -- (285,42) -- (282,42) -- (285,39) -- (283,40) -- (280,40) -- (281,41) -- (279,42) -- (280,40) -- (283,37) -- (281,38) -- (282,36) -- (285,33) -- (291,30) -- (288,30) -- (291,27) -- (297,24) -- (306,24) -- (303,21) -- (309,18) -- (306,18) -- (303,15) -- (303,18) -- (300,18) -- (303,15) -- (309,12) -- (306,12) -- (309,9) -- (307,10) -- (304,10) -- (305,11) -- (303,12) -- (300,12) -- (297,9) -- (297,12) -- (294,12) -- (295,13) -- (295,16) -- (296,14) -- (297,15) -- (295,16) -- (292,16) -- (293,17) -- (291,18) -- (292,16) -- (295,13) -- (293,14) -- (294,12) -- (297,9) -- (303,6) -- (300,6) -- (303,3) -- (301,4) -- (298,4) -- (299,5) -- (297,6) -- (298,4) -- (301,1) -- (299,2) -- (300,0) -- (299,1) -- (297,2) -- (298,2) -- (297,3) -- (295,4) -- (292,4) -- (293,5) -- (291,6) -- (292,6) -- (293,7) -- (293,6) -- (294,6) -- (293,7) -- (291,8) -- (292,8) -- (291,9) -- (289,10) -- (286,10) -- (287,11) -- (285,12) -- (282,12) -- (279,9) -- (279,12) -- (276,12) -- (277,13) -- (277,16) -- (278,14) -- (279,15) -- (277,16) -- (274,16) -- (275,17) -- (273,18) -- (274,18) -- (275,19) -- (275,18) -- (276,18) -- (277,19) -- (277,22) -- (278,20) -- (279,21) -- (279,20) -- (280,18) -- (279,19) -- (279,18) -- (280,18) -- (281,19) -- (281,18) -- (282,18) -- (281,19) -- (279,20) -- (280,20) -- (279,21) -- (277,22) -- (274,22) -- (275,23) -- (273,24) -- (274,24) -- (275,25) -- (275,24) -- (276,24) -- (275,25) -- (273,26) -- (274,26) -- (273,27) -- (273,26) -- (274,24) -- (273,25) -- (273,24) -- (274,22) -- (277,19) -- (275,20) -- (276,18) -- (275,19) -- (273,20) -- (274,20) -- (273,21) -- (273,20) -- (274,18) -- (273,19) -- (273,18) -- (274,16) -- (277,13) -- (275,14) -- (276,12) -- (279,9) -- (285,6) -- (282,6) -- (285,3) -- (283,4) -- (280,4) -- (281,5) -- (279,6) -- (280,4) -- (283,1) -- (281,2) -- (282,0) -- (281,1) -- (279,2) -- (280,2) -- (279,3) -- (277,4) -- (274,4) -- (275,5) -- (273,6) -- (274,6) -- (275,7) -- (275,6) -- (276,6) -- (275,7) -- (273,8) -- (274,8) -- (273,9) -- (273,8) -- (274,6) -- (273,7) -- (273,6) -- (274,4) -- (277,1) -- (275,2) -- (276,0) -- (275,1) -- (273,2) -- (274,2) -- (273,3) -- (273,2) -- (274,0) -- (273,1) -- (273,0);
    	\end{scope}
    \end{tikzpicture}
    }
    \caption{The first four iterations of the three entangled Kochawave curves}
    \label{fig:tri-steps}
\end{figure}
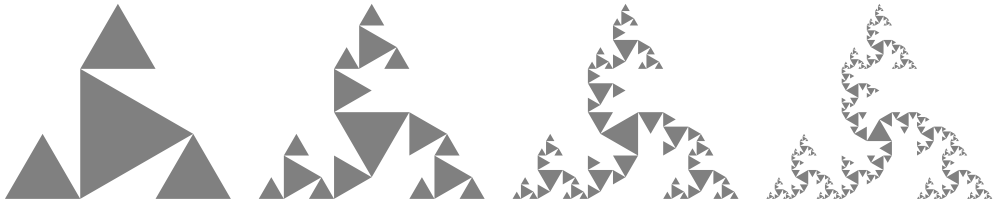

The first four iterations of the construction of the three entangled Kochawave curves are depicted in figure~\ref{fig:tri-steps}.

\begin{figure}[!ht]
  \includegraphics[width=\linewidth]{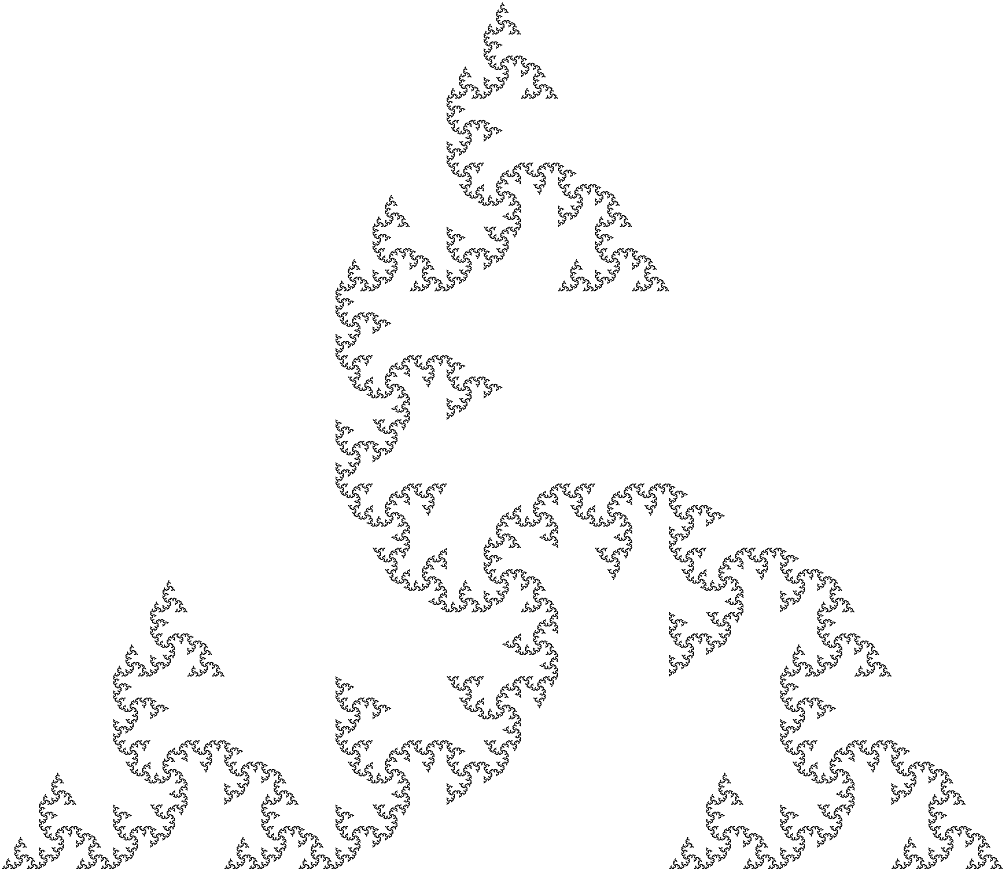}
  \caption{The three entangled Kochawave curves}
  \label{fig:3-kochawaves}
\end{figure}

The three entangled Kochawave curves are the limiting figure as the iterations continue indefinitely (see figure~\ref{fig:3-kochawaves}).
We see that the limiting figure is made up of three
copies of the Kochawave curve we saw in figure~\ref{fig:kochawave}.

\clearpage

\subsection{Construction with an L-system}
\label{cons-l-system}

We can express the Kochawave curve with the following L-system:
\begin{itemize}
  \item alphabet: $\{ S, \ * \ , \ / , \ + \}$,
  \item axiom: $S$,
  \item unique production rule: $S \rightarrow S \ * \ S \ / \ S \ + \ S$.
\end{itemize}

The alphabet is to be interpreted as described in the following table.

\vspace{\baselineskip}

\begin{tabular}{|c|l|}
    \hline 
    \rule[-1ex]{0pt}{2.5ex} Symbol & Interpretation \\ 
    \hline 
    \rule[-1ex]{0pt}{2.5ex} $S$ & Draw one step forward. \\ 
    \rule[-1ex]{0pt}{2.5ex} $*$ & Turn $30^\circ$ left and stretch the next steps by $\sqrt{3}$. \\ 
    \rule[-1ex]{0pt}{2.5ex} $/$ & Turn $150^\circ$ right and stretch the next steps by $1/\sqrt{3}$. \\ 
    \rule[-1ex]{0pt}{2.5ex} $+$ & Turn $120^\circ$ left. \\ 
    \hline 
\end{tabular} 

\vspace{\baselineskip}

The length of the $k^{th}$ step is $\sqrt{3}^{u_{k-1}}$,
where $u_{k-1}$ denotes the number of $1$'s in the base-4 representation of $k-1$.

\subsection{Numerical construction}
\label{cons-numeric}

Let:
\begin{equation*}
    z_0 = 0, \ 
    z_{k+1} = z_k + \frac{(1+\omega)^{u_k}}{\omega^{2v_k}}
\end{equation*}
where $\omega = e^{\frac{i\pi}{3}}$, and $u_k$ and $v_k$ denote the number of $1$'s and of $2$'s, respectively,
in the base-4 expansion of $k$.

The coordinates of the sequence of points $\{z_n\}_{n \geq 0}$
on a hexagonal lattice are given by the OEIS~\cite{oeis} sequences \href{https://oeis.org/A335380}{A335380}
and \href{https://oeis.org/A335381}{A335381}.

The set of endpoints of the segments after $n$ iterations
of the construction described in section~\ref{cons-segments}
is given by:
\begin{equation*}
    K_n = \left\{ \frac{z_k}{3^n}, \ k = 0 \twodots 4^n \right\}
\end{equation*} 

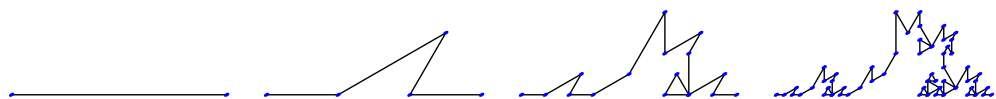
\begin{figure}[!ht]
    \resizebox{13cm}{!}
    {
    \begin{tikzpicture}
    	\begin{scope}[yscale=.87,xslant=.5]
    		\draw[line width=5pt] (0,0) -- (27,0);
    		\draw[line width=5pt] (32,0) -- (41,0) -- (50,9) -- (50,0) -- (59,0);
    		\draw[line width=5pt] (64,0) -- (67,0) -- (70,3) -- (70,0) -- (73,0) -- (76,3) -- (76,12) -- (79,6) -- (82,9) -- (82,6) -- (85,0) -- (82,3) -- (82,0) -- (85,0) -- (88,3) -- (88,0) -- (91,0);
    		\draw[line width=5pt] (96,0) -- (97,0) -- (98,1) -- (98,0) -- (99,0) -- (100,1) -- (100,4) -- (101,2) -- (102,3) -- (102,2) -- (103,0) -- (102,1) -- (102,0) -- (103,0) -- (104,1) -- (104,0) -- (105,0) -- (106,1) -- (106,4) -- (107,2) -- (108,3) -- (108,6) -- (105,12) -- (108,9) -- (108,12) -- (109,10) -- (112,7) -- (110,8) -- (111,6) -- (112,7) -- (112,10) -- (113,8) -- (114,9) -- (114,8) -- (115,6) -- (114,7) -- (114,6) -- (115,4) -- (118,1) -- (116,2) -- (117,0) -- (116,1) -- (114,2) -- (115,2) -- (114,3) -- (114,2) -- (115,0) -- (114,1) -- (114,0) -- (115,0) -- (116,1) -- (116,0) -- (117,0) -- (118,1) -- (118,4) -- (119,2) -- (120,3) -- (120,2) -- (121,0) -- (120,1) -- (120,0) -- (121,0) -- (122,1) -- (122,0) -- (123,0);
    		\begin{scope}[fill=blue]
    			\fill (0,0) circle (3mm); \fill (27,0) circle (3mm);
    			\fill (32,0) circle (3mm); \fill (41,0) circle (3mm); \fill (50,9) circle (3mm); \fill (50,0) circle (3mm); \fill (59,0) circle (3mm);
    			\fill (64,0) circle (3mm); \fill (67,0) circle (3mm); \fill (70,3) circle (3mm); \fill (70,0) circle (3mm); \fill (73,0) circle (3mm); \fill (76,3) circle (3mm); \fill (76,12) circle (3mm); \fill (79,6) circle (3mm); \fill (82,9) circle (3mm); \fill (82,6) circle (3mm); \fill (85,0) circle (3mm); \fill (82,3) circle (3mm); \fill (82,0) circle (3mm); \fill (85,0) circle (3mm); \fill (88,3) circle (3mm); \fill (88,0) circle (3mm); \fill (91,0) circle (3mm);
    			\fill (96,0) circle (3mm); \fill (97,0) circle (3mm); \fill (98,1) circle (3mm); \fill (98,0) circle (3mm); \fill (99,0) circle (3mm); \fill (100,1) circle (3mm); \fill (100,4) circle (3mm); \fill (101,2) circle (3mm); \fill (102,3) circle (3mm); \fill (102,2) circle (3mm); \fill (103,0) circle (3mm); \fill (102,1) circle (3mm); \fill (102,0) circle (3mm); \fill (103,0) circle (3mm); \fill (104,1) circle (3mm); \fill (104,0) circle (3mm); \fill (105,0) circle (3mm); \fill (106,1) circle (3mm); \fill (106,4) circle (3mm); \fill (107,2) circle (3mm); \fill (108,3) circle (3mm); \fill (108,6) circle (3mm); \fill (105,12) circle (3mm); \fill (108,9) circle (3mm); \fill (108,12) circle (3mm); \fill (109,10) circle (3mm); \fill (112,7) circle (3mm); \fill (110,8) circle (3mm); \fill (111,6) circle (3mm); \fill (112,7) circle (3mm); \fill (112,10) circle (3mm); \fill (113,8) circle (3mm); \fill (114,9) circle (3mm); \fill (114,8) circle (3mm); \fill (115,6) circle (3mm); \fill (114,7) circle (3mm); \fill (114,6) circle (3mm); \fill (115,4) circle (3mm); \fill (118,1) circle (3mm); \fill (116,2) circle (3mm); \fill (117,0) circle (3mm); \fill (116,1) circle (3mm); \fill (114,2) circle (3mm); \fill (115,2) circle (3mm); \fill (114,3) circle (3mm); \fill (114,2) circle (3mm); \fill (115,0) circle (3mm); \fill (114,1) circle (3mm); \fill (114,0) circle (3mm); \fill (115,0) circle (3mm); \fill (116,1) circle (3mm); \fill (116,0) circle (3mm); \fill (117,0) circle (3mm); \fill (118,1) circle (3mm); \fill (118,4) circle (3mm); \fill (119,2) circle (3mm); \fill (120,3) circle (3mm); \fill (120,2) circle (3mm); \fill (121,0) circle (3mm); \fill (120,1) circle (3mm); \fill (120,0) circle (3mm); \fill (121,0) circle (3mm); \fill (122,1) circle (3mm); \fill (122,0) circle (3mm); \fill (123,0) circle (3mm);
    		\end{scope}
    	\end{scope}
    \end{tikzpicture}
    }
    \caption{The sets of points $K_0$ to $K_3$ (in blue)}
    \label{fig:points}
\end{figure}

The first few sets are depicted in figure~\ref{fig:points}.

As $z_{4n} = 3 z_n$, $K_n \subset K_{n+1}$.

The closure of the limit $\lim_{n \to \infty} K_n$ is the Kochawave curve.

\section{Properties}

\subsection{Length of the Kochawave curve}
After $n$ iterations of the construction described in section~\ref{cons-segments},
the length $L_n$ of the structure is given by:
\begin{equation*}
    L_n = \left(1 + \frac{1}{\sqrt{3}} \right)^n
\end{equation*}
Note that $L_0 = 1$ as we started with a unit segment.

As $\lim_{n \to \infty} L_n = \infty$, the Kochawave curve has infinite length.

\subsection{Area of the three entangled Kochawave curves}
After $n$ iterations of the construction described in section~\ref{cons-triangles},
the area $T_n$ of the structure is given by:
\begin{equation*}
    T_n = \frac{\sqrt{3}}{4} \cdot \left(\frac{2}{3}\right)^n
\end{equation*}
Note that $T_0 = \frac{\sqrt{3}}{4}$ as we started with a unit sided equilateral triangle.

As $\lim_{n \to \infty} T_n = 0$, the three entangled Kochawave curves have empty area.
By contrast, the Koch antisnowflake has positive area.

\subsection{Area of the Kochawave curve}
\label{prop-area}
After $n$ iterations of the construction described in section~\ref{cons-segments},
the area $A_n$ enclosed between the structure and the initial segment satisfies:
\begin{equation*}
    3A_n + T_n = \frac{\sqrt{3}}{4}
\end{equation*}

Hence:
\begin{equation*}
    A_n = \frac{\frac{\sqrt{3}}{4} - T_n}{3}
        = \frac{1}{4\sqrt{3}} \cdot \left[1 - \left(\frac{2}{3}\right)^n \right]
\end{equation*}

The area $A$ of the Kochawave curve is given by:
\begin{equation*}
    A = lim_{n \to \infty} A_n = \frac{1}{4\sqrt{3}}
\end{equation*}

\begin{figure}[!ht]
    \includegraphics[width=\linewidth]{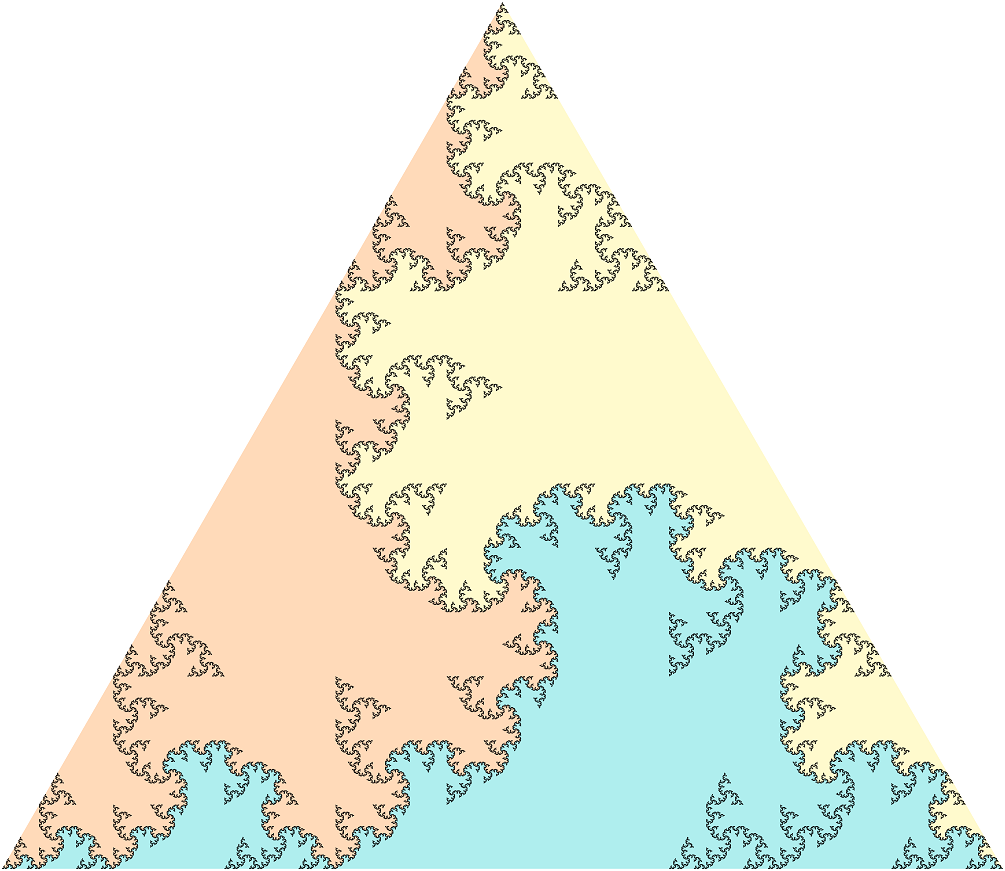}
    \caption{An equilateral triangle divided into three Kochawave curves}
    \label{fig:area}
\end{figure}

This is one third of the area of a unit sided equilateral triangle,
and indeed, we can divide an equilateral triangle
into three Kochawave curves (see figure~\ref{fig:area}).

\subsection{Height of the Kochawave curve}
After $n$ iterations of the construction described in section~\ref{cons-segments},
the height $H_n$ of the structure is given by:
\begin{equation*}
    H_n = \max_{p \in K_n} Im(p)
        = \left\{ \begin{array}{cl}
                0                      & \mathrm{if} \ n = 0 \\
                \frac{1}{2\sqrt{3}}    & \mathrm{if} \ n = 1 \\
                \frac{2}{3\sqrt{3}}    & \mathrm{if} \ n > 1
          \end{array} \right.
\end{equation*}

The height $H$ of the Kochawave curve is given by:
\begin{equation*}
    H = lim_{n \to \infty} H_n = \frac{2}{3\sqrt{3}}
\end{equation*}

The leftmost  point of height $H$ has position $\frac{4 + 4\omega}{9}$.

The rightmost point of height $H$ has position $\frac{5 + 4\omega}{9}$.

\subsection{Centroid of the Kochawave curve}

\begin{figure}[!ht]
  \includegraphics[width=\linewidth]{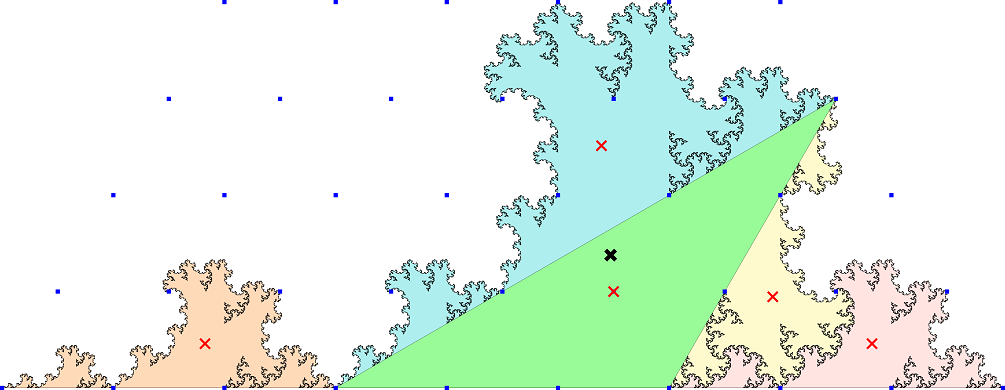}
  \caption{The dissection of a Kochawave curve into five parts
          (red marks correspond to centroids of parts, the black mark corresponds to the centroid $m$)}
  \label{fig:dissection}
\end{figure}

We can dissect a Kochawave curve, say of weight $1$, into five parts:
\begin{itemize}
    \item    one copy                  of the curve of weight $\frac{1}{3}$,
    \item    three copies of the curve of weight $\frac{1}{9}$,
    \item    one triangle              of weight $\frac{1}{3}$.
\end{itemize}
as depicted in figure~\ref{fig:dissection}
(this dissection derives directly from the construction described in section~\ref{cons-segments}).

This leads to the following linear equation for the centroid $m$
of the Kochawave curve:
\begin{equation*}
m = \frac{1}{9} \cdot \frac{                                           m}{3}
  + \frac{1}{3} \cdot \frac{1                           + (1 + \omega) m}{3}
  + \frac{1}{9} \cdot \frac{2 + \omega                  - \omega       m}{3}
  + \frac{1}{9} \cdot \frac{2                           +              m}{3}
  + \frac{1}{3} \cdot \frac{\frac{5}{3} + \frac{\omega}{3}              }{3}
\end{equation*}

Hence the centroid $m$ of the Kochawave curve is:
\begin{equation*}
    m = \frac{59+17\omega}{111}
\end{equation*}

\subsection{Volume of solid of revolution}

The volume $V$ of the solid of revolution of the Kochawave curve
around its base unit segment satisfies:

\begin{equation*}
    V = 2 \pi \cdot Im(m) \cdot A = \frac{17}{444} \pi
\end{equation*}

\subsection{Hausdorff dimension of the Kochawave curve}

The Kochawave curve satisfies the open set condition;
its Hausdorff dimension $d$ is the unique real solution of the following equation:
\begin{equation*}
		\frac{1}{3^{d-1}} + \frac{1}{\sqrt{3}^d} = 1
\end{equation*}
and has exact value:
\begin{equation*}
		d = 2 \cdot \frac{\log \frac{1+\sqrt{13}}{2}}{\log 3}
\end{equation*}
which is approximately $1.518$.

\subsection{Symmetries}

Contrary to the Koch curve, the Kochawave curve has no particular symmetry (see figure~\ref{fig:kochawave}).
As for the Koch antisnowflake, the three entangled Kochawave curves have 3-fold rotational symmetry (see figure~\ref{fig:3-kochawaves}).

\subsection{Cantor set underlying the Kochawave curve}

When building the Kochawave curve starting from one unit segment
as described in section~\ref{cons-segments},
the remainder of this initial segment forms a ternary Cantor set.
This property is shared with the Koch curve.

\begin{figure}[!ht]
    \resizebox{13cm}{!}
    {
    \begin{tikzpicture}
    	\begin{scope}[yscale=.87,xslant=.5]
    		\draw[line width=5pt] (0,0) -- (27,0) -- (54,27) -- (54,0) -- (81,0);
    		\draw[line width=5pt] (91,0) -- (100,0) -- (109,9) -- (109,0) -- (118,0) -- (127,9) -- (127,36) -- (136,18) -- (145,27) -- (145,18) -- (154,0) -- (145,9) -- (145,0) -- (154,0) -- (163,9) -- (163,0) -- (172,0);
    		\draw[line width=5pt] (182,0) -- (185,0) -- (188,3) -- (188,0) -- (191,0) -- (194,3) -- (194,12) -- (197,6) -- (200,9) -- (200,6) -- (203,0) -- (200,3) -- (200,0) -- (203,0) -- (206,3) -- (206,0) -- (209,0) -- (212,3) -- (212,12) -- (215,6) -- (218,9) -- (218,18) -- (209,36) -- (218,27) -- (218,36) -- (221,30) -- (230,21) -- (224,24) -- (227,18) -- (230,21) -- (230,30) -- (233,24) -- (236,27) -- (236,24) -- (239,18) -- (236,21) -- (236,18) -- (239,12) -- (248,3) -- (242,6) -- (245,0) -- (242,3) -- (236,6) -- (239,6) -- (236,9) -- (236,6) -- (239,0) -- (236,3) -- (236,0) -- (239,0) -- (242,3) -- (242,0) -- (245,0) -- (248,3) -- (248,12) -- (251,6) -- (254,9) -- (254,6) -- (257,0) -- (254,3) -- (254,0) -- (257,0) -- (260,3) -- (260,0) -- (263,0);
    		\draw[line width=5pt] (273,0) -- (274,0) -- (275,1) -- (275,0) -- (276,0) -- (277,1) -- (277,4) -- (278,2) -- (279,3) -- (279,2) -- (280,0) -- (279,1) -- (279,0) -- (280,0) -- (281,1) -- (281,0) -- (282,0) -- (283,1) -- (283,4) -- (284,2) -- (285,3) -- (285,6) -- (282,12) -- (285,9) -- (285,12) -- (286,10) -- (289,7) -- (287,8) -- (288,6) -- (289,7) -- (289,10) -- (290,8) -- (291,9) -- (291,8) -- (292,6) -- (291,7) -- (291,6) -- (292,4) -- (295,1) -- (293,2) -- (294,0) -- (293,1) -- (291,2) -- (292,2) -- (291,3) -- (291,2) -- (292,0) -- (291,1) -- (291,0) -- (292,0) -- (293,1) -- (293,0) -- (294,0) -- (295,1) -- (295,4) -- (296,2) -- (297,3) -- (297,2) -- (298,0) -- (297,1) -- (297,0) -- (298,0) -- (299,1) -- (299,0) -- (300,0) -- (301,1) -- (301,4) -- (302,2) -- (303,3) -- (303,6) -- (300,12) -- (303,9) -- (303,12) -- (304,10) -- (307,7) -- (305,8) -- (306,6) -- (307,7) -- (307,10) -- (308,8) -- (309,9) -- (309,12) -- (306,18) -- (309,15) -- (309,18) -- (306,24) -- (297,33) -- (303,30) -- (300,36) -- (303,33) -- (309,30) -- (306,30) -- (309,27) -- (309,30) -- (306,36) -- (309,33) -- (309,36) -- (310,34) -- (313,31) -- (311,32) -- (312,30) -- (315,27) -- (321,24) -- (318,24) -- (321,21) -- (319,22) -- (316,22) -- (317,23) -- (315,24) -- (316,22) -- (319,19) -- (317,20) -- (318,18) -- (319,19) -- (319,22) -- (320,20) -- (321,21) -- (321,24) -- (318,30) -- (321,27) -- (321,30) -- (322,28) -- (325,25) -- (323,26) -- (324,24) -- (325,25) -- (325,28) -- (326,26) -- (327,27) -- (327,26) -- (328,24) -- (327,25) -- (327,24) -- (328,22) -- (331,19) -- (329,20) -- (330,18) -- (329,19) -- (327,20) -- (328,20) -- (327,21) -- (327,20) -- (328,18) -- (327,19) -- (327,18) -- (328,16) -- (331,13) -- (329,14) -- (330,12) -- (333,9) -- (339,6) -- (336,6) -- (339,3) -- (337,4) -- (334,4) -- (335,5) -- (333,6) -- (334,4) -- (337,1) -- (335,2) -- (336,0) -- (335,1) -- (333,2) -- (334,2) -- (333,3) -- (331,4) -- (328,4) -- (329,5) -- (327,6) -- (328,6) -- (329,7) -- (329,6) -- (330,6) -- (329,7) -- (327,8) -- (328,8) -- (327,9) -- (327,8) -- (328,6) -- (327,7) -- (327,6) -- (328,4) -- (331,1) -- (329,2) -- (330,0) -- (329,1) -- (327,2) -- (328,2) -- (327,3) -- (327,2) -- (328,0) -- (327,1) -- (327,0) -- (328,0) -- (329,1) -- (329,0) -- (330,0) -- (331,1) -- (331,4) -- (332,2) -- (333,3) -- (333,2) -- (334,0) -- (333,1) -- (333,0) -- (334,0) -- (335,1) -- (335,0) -- (336,0) -- (337,1) -- (337,4) -- (338,2) -- (339,3) -- (339,6) -- (336,12) -- (339,9) -- (339,12) -- (340,10) -- (343,7) -- (341,8) -- (342,6) -- (343,7) -- (343,10) -- (344,8) -- (345,9) -- (345,8) -- (346,6) -- (345,7) -- (345,6) -- (346,4) -- (349,1) -- (347,2) -- (348,0) -- (347,1) -- (345,2) -- (346,2) -- (345,3) -- (345,2) -- (346,0) -- (345,1) -- (345,0) -- (346,0) -- (347,1) -- (347,0) -- (348,0) -- (349,1) -- (349,4) -- (350,2) -- (351,3) -- (351,2) -- (352,0) -- (351,1) -- (351,0) -- (352,0) -- (353,1) -- (353,0) -- (354,0);
    		\begin{scope}[line width=1cm,blue]
    			\draw (0,0) -- (27,0); \draw (54,0) -- (81,0);
    			\draw (91,0) -- (100,0); \draw (109,0) -- (118,0); \draw (145,0) -- (154,0); \draw (163,0) -- (172,0);
    			\draw (182,0) -- (185,0); \draw (188,0) -- (191,0); \draw (200,0) -- (203,0); \draw (206,0) -- (209,0); \draw (236,0) -- (239,0); \draw (242,0) -- (245,0); \draw (254,0) -- (257,0); \draw (260,0) -- (263,0);
    			\draw (273,0) -- (274,0); \draw (275,0) -- (276,0); \draw (279,0) -- (280,0); \draw (281,0) -- (282,0); \draw (291,0) -- (292,0); \draw (293,0) -- (294,0); \draw (297,0) -- (298,0); \draw (299,0) -- (300,0); \draw (327,0) -- (328,0); \draw (329,0) -- (330,0); \draw (333,0) -- (334,0); \draw (335,0) -- (336,0); \draw (345,0) -- (346,0); \draw (347,0) -- (348,0); \draw (351,0) -- (352,0); \draw (353,0) -- (354,0);
    		\end{scope}
    	\end{scope}
    \end{tikzpicture}
    }
    \caption{Remainder of the initial segment after one to four iterations (in blue)}
    \label{fig:cantor}
\end{figure}
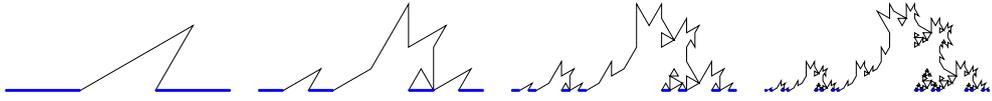

\subsection{Connectivity}

The Kochawave curve and the three entangled Kochawave curves are simply connected.

\subsection{Left boundary}

\begin{figure}[!ht]
    \resizebox{13cm}{!}
    {
    \begin{tikzpicture}
    	\begin{scope}[yscale=.87,xslant=.5]
    		\draw (0,0) -- (1,0) -- (2,1) -- (2,0) -- (3,0) -- (4,1) -- (4,4) -- (5,2) -- (6,3) -- (6,2) -- (7,0) -- (8,1) -- (8,0) -- (9,0) -- (10,1) -- (10,4) -- (11,2) -- (12,3) -- (12,6) -- (9,12) -- (12,9) -- (12,12) -- (13,10) -- (16,7) -- (16,10) -- (17,8) -- (18,9) -- (18,8) -- (19,6) -- (18,7) -- (18,6) -- (19,4) -- (22,1) -- (22,4) -- (23,2) -- (24,3) -- (24,2) -- (25,0) -- (26,1) -- (26,0) -- (27,0) -- (28,1) -- (28,4) -- (29,2) -- (30,3) -- (30,6) -- (27,12) -- (30,9) -- (30,12) -- (31,10) -- (34,7) -- (34,10) -- (35,8) -- (36,9) -- (36,12) -- (33,18) -- (36,15) -- (36,18) -- (33,24) -- (24,33) -- (30,30) -- (27,36) -- (30,33) -- (36,30) -- (33,36) -- (36,33) -- (36,36) -- (37,34) -- (40,31) -- (38,32) -- (39,30) -- (42,27) -- (48,24) -- (45,30) -- (48,27) -- (48,30) -- (49,28) -- (52,25) -- (52,28) -- (53,26) -- (54,27) -- (54,26) -- (55,24) -- (54,25) -- (54,24) -- (55,22) -- (58,19) -- (56,20) -- (57,18) -- (56,19) -- (54,20) -- (55,18) -- (54,19) -- (54,18) -- (55,16) -- (58,13) -- (56,14) -- (57,12) -- (60,9) -- (66,6) -- (63,12) -- (66,9) -- (66,12) -- (67,10) -- (70,7) -- (70,10) -- (71,8) -- (72,9) -- (72,8) -- (73,6) -- (72,7) -- (72,6) -- (73,4) -- (76,1) -- (76,4) -- (77,2) -- (78,3) -- (78,2) -- (79,0) -- (80,1) -- (80,0) -- (81,0);
    		\begin{scope}[line width=2mm,green]
    			\draw (7,0) -- (6,1) -- (6,0) -- cycle;
    			\draw (16,7) -- (14,8) -- (15,6) -- cycle;
    			\draw (22,1) -- (20,2) -- (21,0) -- (20,1) -- (18,2) -- (19,2) -- (18,3) -- (18,2) -- (19,0) -- (18,1) -- (18,0) -- (19,0) -- (20,1) -- (20,0) -- (21,0) -- cycle;
    			\draw (25,0) -- (24,1) -- (24,0) -- cycle;
    			\draw (34,7) -- (32,8) -- (33,6) -- cycle;
    			\draw (36,30) -- (33,30) -- (36,27) -- cycle;
    			\draw (48,24) -- (45,24) -- (48,21) -- (46,22) -- (43,22) -- (44,23) -- (42,24) -- (43,22) -- (46,19) -- (44,20) -- (45,18) -- (46,19) -- (46,22) -- (47,20) -- (48,21) -- cycle;
    			\draw (52,25) -- (50,26) -- (51,24) -- cycle;
    			\draw (54,20) -- (55,20) -- (54,21) -- cycle;
    			\draw (66,6) -- (63,6) -- (66,3) -- (64,4) -- (61,4) -- (62,5) -- (60,6) -- (61,4) -- (64,1) -- (62,2) -- (63,0) -- (62,1) -- (60,2) -- (61,2) -- (60,3) -- (58,4) -- (55,4) -- (56,5) -- (54,6) -- (55,6) -- (56,7) -- (56,6) -- (57,6) -- (56,7) -- (54,8) -- (55,8) -- (54,9) -- (54,8) -- (55,6) -- (54,7) -- (54,6) -- (55,4) -- (58,1) -- (56,2) -- (57,0) -- (56,1) -- (54,2) -- (55,2) -- (54,3) -- (54,2) -- (55,0) -- (54,1) -- (54,0) -- (55,0) -- (56,1) -- (56,0) -- (57,0) -- (58,1) -- (58,4) -- (59,2) -- (60,3) -- (60,2) -- (61,0) -- (60,1) -- (60,0) -- (61,0) -- (62,1) -- (62,0) -- (63,0) -- (64,1) -- (64,4) -- (65,2) -- (66,3) -- cycle;
    			\draw (70,7) -- (68,8) -- (69,6) -- cycle;
    			\draw (76,1) -- (74,2) -- (75,0) -- (74,1) -- (72,2) -- (73,2) -- (72,3) -- (72,2) -- (73,0) -- (72,1) -- (72,0) -- (73,0) -- (74,1) -- (74,0) -- (75,0) -- cycle;
    			\draw (79,0) -- (78,1) -- (78,0) -- cycle;
    		\end{scope}
    	\end{scope}
    \end{tikzpicture}
    }
    \caption{The closed loops after four iterations of the construction of the Kochawave curve (in green)}
    \label{fig:closed-loops}
\end{figure}
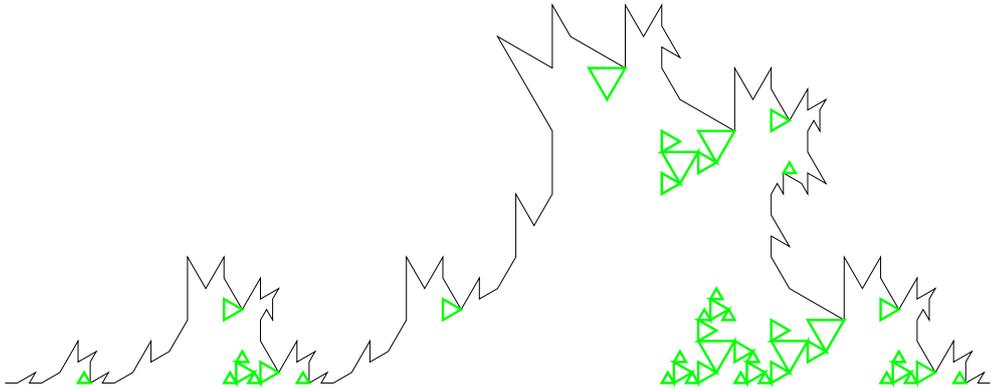

During the construction of the Kochawave curve seen in section~\ref{cons-segments},
we notice closed loops made of equilateral triangles.
Each of these equilateral triangles will give rise
to three small entangled Kochawave curves,
and participate in the fractal nature of the Kochawave curve (see figure~\ref{fig:closed-loops}).

\begin{figure}[!ht]
  \includegraphics[width=\linewidth]{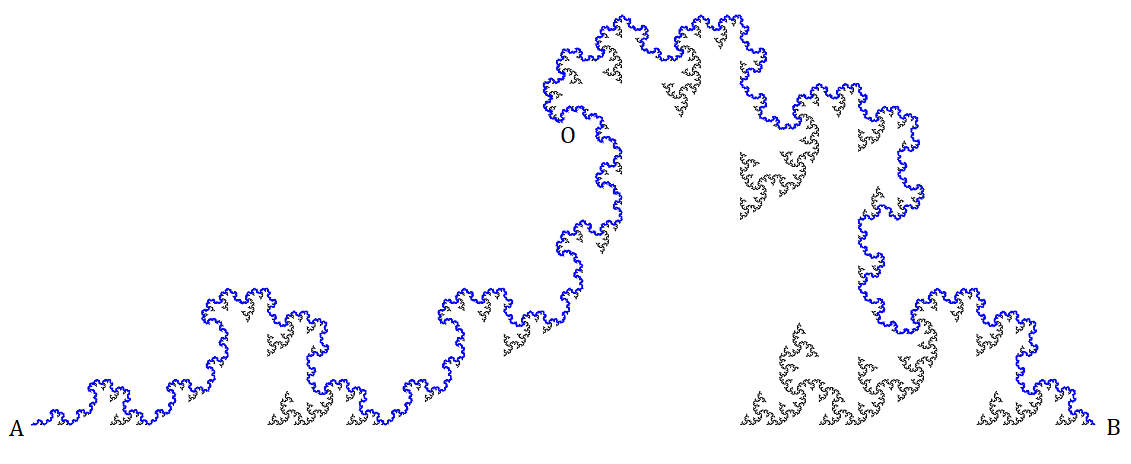}
  \caption{The curve $\mathcal{C}$ (in blue)}
  \label{fig:curve-c}
\end{figure}

If we remove these closed loops, we obtain a simple curve $\mathcal{C}$,
with the same area (see section~\ref{prop-area}) as the corresponding Kochawave curve.
The curve $\mathcal{C}$ can be split into two parts,
$AO$ and $BO$ (where $O$ is the center of the equilateral triangle with base $AB$);
$AO$, when rotated $120^{\circ}$ counterclockwise around $O$, gives $BO$ (see figure~\ref{fig:curve-c}).

\begin{figure}[!ht]
    \resizebox{13cm}{!}
    {
    \begin{tikzpicture}
    	\begin{scope}[yscale=.87,xslant=.5]
    		\begin{scope}[thick]
    			\draw (0,0) node[below] {\tiny A} -- (3,0) node[below] {\tiny D};

    			\draw[->,gray] (3.5,0.5) -> (4.5,0.5);

    			\draw[fill,gray!20] (7,0) -- (8,0) -- (7,1);
    			\draw (5,0) node[below] {\tiny A} -- (6,0) node[below] {\tiny B} -- (7,1) node[above] {\tiny C};
    			\draw[blue,text=black] (7,1) -- (7,0);
    			\draw[red,text=black] (7,0) -- (8,0) node[below] {\tiny D};

    			\draw[fill,gray!20] (2,-1) -- (2,-4) -- (5,-4);
    			\draw[blue,text=black] (2,-1) node[above] {\tiny E} -- (2,-4);
    			\draw[red,text=black] (2,-4) -- (5,-4) node[below] {\tiny H};

    			\draw[->,gray] (5.25,-3) -> (6.25,-3);

    			\draw[fill,gray!20] (7,-2) -- (9,-3) -- (8,-4);
    			\draw[fill,gray!20] (9,-3) -- (9,-4) -- (10,-4);

    			\draw (7,-1) node[above] {\tiny E} -- (7,-2) node[above left] {\tiny F};
    			\draw[blue] (7,-2) -- (8,-4);
    			\draw[red,text=black] (8,-4) -- (9,-3) node[above] {\tiny G};
    			\draw[blue] (9,-3) -- (9,-4);
    			\draw[red,text=black] (9,-4) -- (10,-4) node[below] {\tiny H};
    		\end{scope}
    	\end{scope}
    \end{tikzpicture}
    }
    \caption{The construction rules for the curve $\mathcal{C}$}
    \label{fig:c-one-step}
\end{figure}
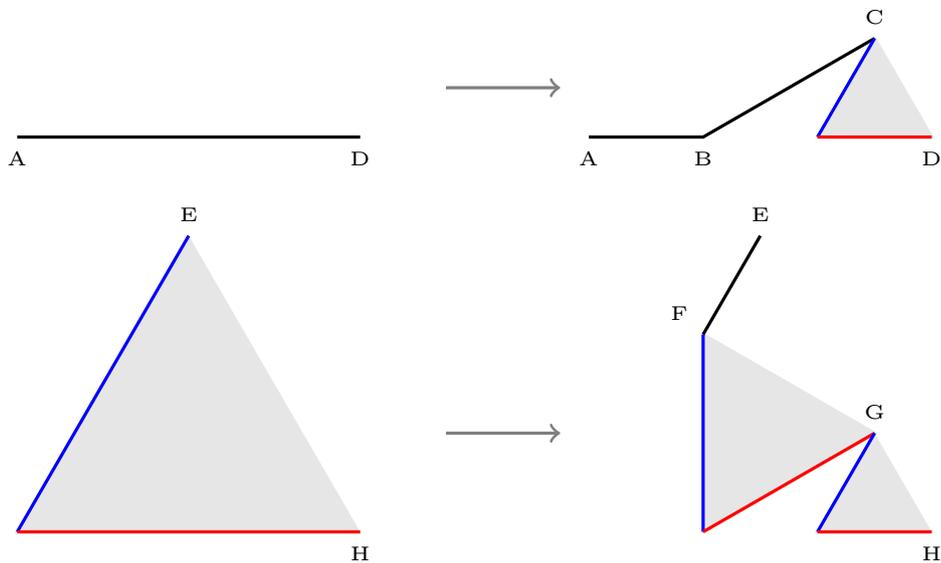

We can construct the curve $\mathcal{C}$ starting with a single black unit segment,
and applying the substitution rules depicted in figure~\ref{fig:c-one-step}.

After $n$ iterations, the total length of the construction, say $c_n$, satisfies:
\begin{equation*}
    \left(\frac{1}{2} + \frac{1}{\sqrt{3}}\right)^n \leq c_n \leq \left(1 + \frac{1}{\sqrt{3}}\right)^n
\end{equation*}
These inequalities are deduced from the two construction rules in figure~\ref{fig:c-one-step}.

As $\frac{1}{2} + \frac{1}{\sqrt{3}} > 1$, $\lim_{n \to \infty} c_n = \infty$, and the curve $\mathcal{C}$ has infinite length.

\newpage

\section{Tessellations of the plane}

\subsection{Tiles}
\label{tess-tile}

We describe five tiles whose edges are made of two, three or four Kochawave curves.
These tiles can cover the Euclidian plane (see section~\ref{tess-plane}).

\subsubsection{Biface antisymmetrical tile}

\begin{figure}[!ht]
  \begin{center}
  \includegraphics[height=3cm, keepaspectratio]{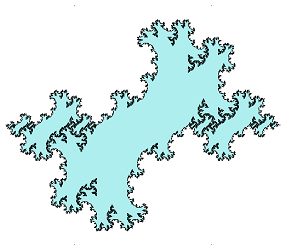}
  \end{center}
  \caption{The biface antisymmetrical tile}
  \label{fig:tile-antisym}
\end{figure}

The biface antisymmetrical tile is made of two Kochawave curves
arranged head to tail as depicted in figure~\ref{fig:tile-antisym}.

This tile can be obtained by arranging $2^{k-1}$ copies of
the rhomboid tile (see figure~\ref{fig:tile-rhomb})
scaled by a factor $1/3^k$ for $k \ge 1$.

\subsubsection{Biface symmetrical tile}

\begin{figure}[!ht]
  \begin{center}
  \includegraphics[height=3cm, keepaspectratio]{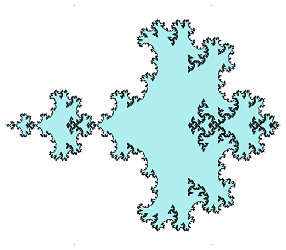}
  \end{center}
  \caption{The biface symmetrical tile}
  \label{fig:tile-sym}
\end{figure}

The biface symmetrical tile is made of two Kochawave curves
arranged symmetrically as depicted in figure~\ref{fig:tile-sym}.

This tile can be obtained by arranging $2^{k-1}$ copies of
the dart tile (see figure~\ref{fig:tile-dart}) scaled by a factor $1/3^k$ for $k \ge 1$.

\subsubsection{Triangular tile}

\begin{figure}[!ht]
  \begin{center}
  \includegraphics[height=3cm, keepaspectratio]{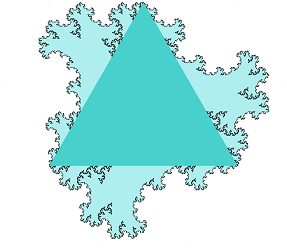}
  \end{center}
  \caption{The triangular tile}
  \label{fig:tile-tri}
\end{figure}

The triangular tile is made of three Kochawave curves
arranged around an equilateral triangle as depicted in figure~\ref{fig:tile-tri}.

\subsubsection{Rhomboidal tile}

\begin{figure}[!ht]
  \begin{center}
  \includegraphics[height=3cm, keepaspectratio]{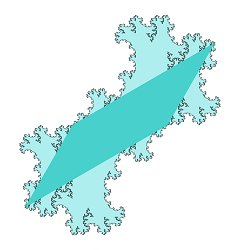}
  \end{center}
  \caption{The rhomboidal tile}
  \label{fig:tile-rhomb}
\end{figure}

The rhomboidal tile is made of four Kochawave curves
arranged around a rhomboid as depicted in figure~\ref{fig:tile-rhomb}.

\subsubsection{Dart tile}

\begin{figure}[!ht]
  \begin{center}
  \includegraphics[height=3cm, keepaspectratio]{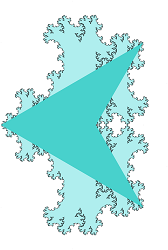}
  \end{center}
  \caption{The dart tile}
  \label{fig:tile-dart}
\end{figure}

The dart tile is made of four Kochawave curves
arranged around a dart as depicted in figure~\ref{fig:tile-dart}.

\subsection{Tessellations}
\label{tess-plane}

We describe five ways to cover the Euclidian plane with the tiles seen in section~\ref{tess-tile}.

\subsubsection{With biface antisymmetrical tiles}

\begin{figure}[!ht]
  \begin{center}
  \includegraphics[height=13cm, keepaspectratio]{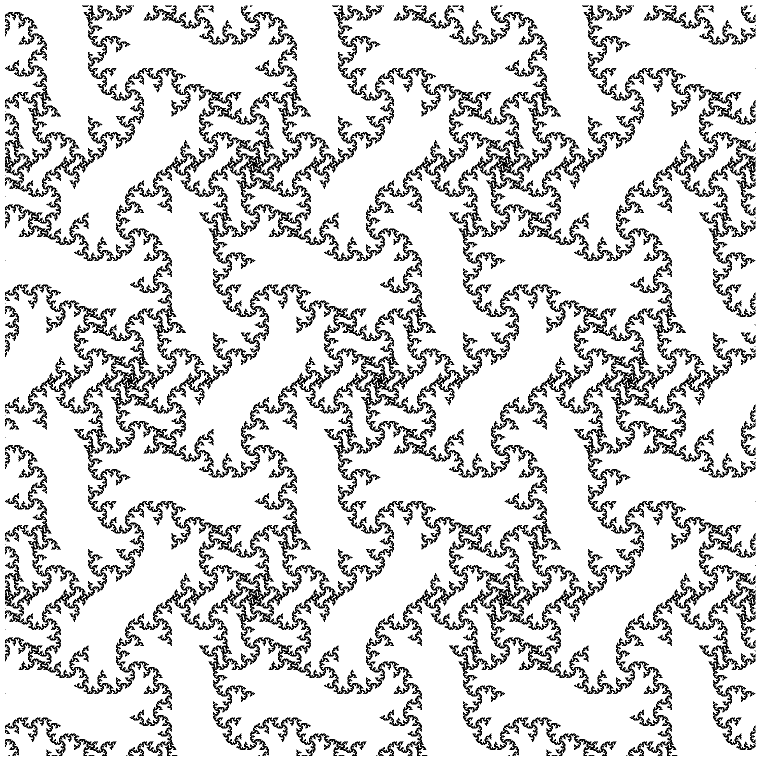}
  \end{center}
  \caption{Covering the plane with biface antisymmetrical tiles}
  \label{fig:plane-antisym}
\end{figure}

We can cover the plane with copies of the biface antisymmetrical tile in one size (see figure~\ref{fig:plane-antisym}).
This covering is periodic.

\clearpage

\subsubsection{With biface antisymmetrical tiles}

\begin{figure}[!ht]
  \begin{center}
  \includegraphics[height=13cm, keepaspectratio]{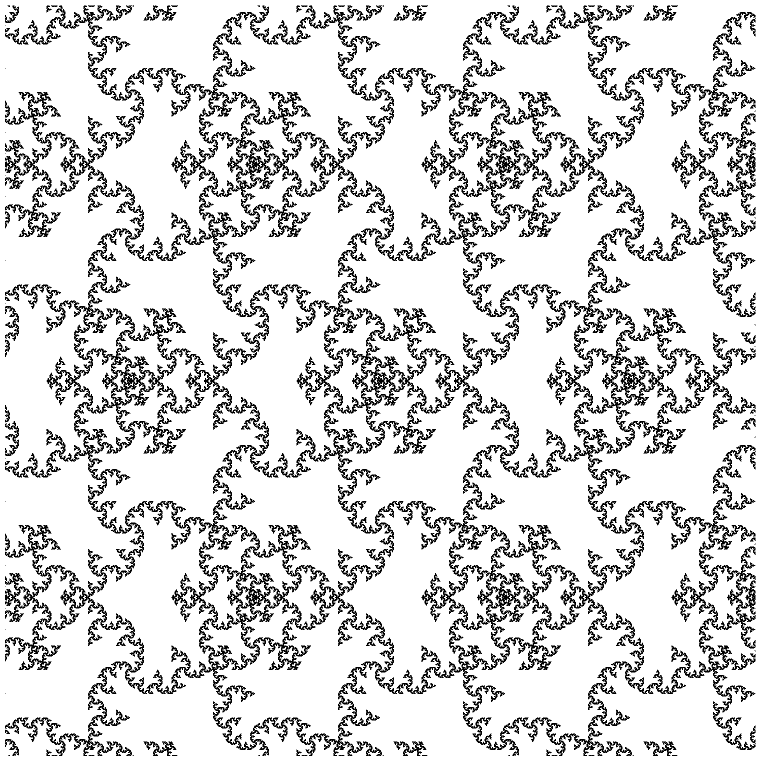}
  \end{center}
  \caption{Covering the plane with biface symmetrical tiles}
  \label{fig:plane-sym}
\end{figure}

We can cover the plane with copies of the biface antisymmetrical tile in one size (see figure~\ref{fig:plane-sym}).
This covering is periodic.

\clearpage

\subsubsection{With triangular tiles}

\begin{figure}[!ht]
  \begin{center}
  \includegraphics[height=13cm, keepaspectratio]{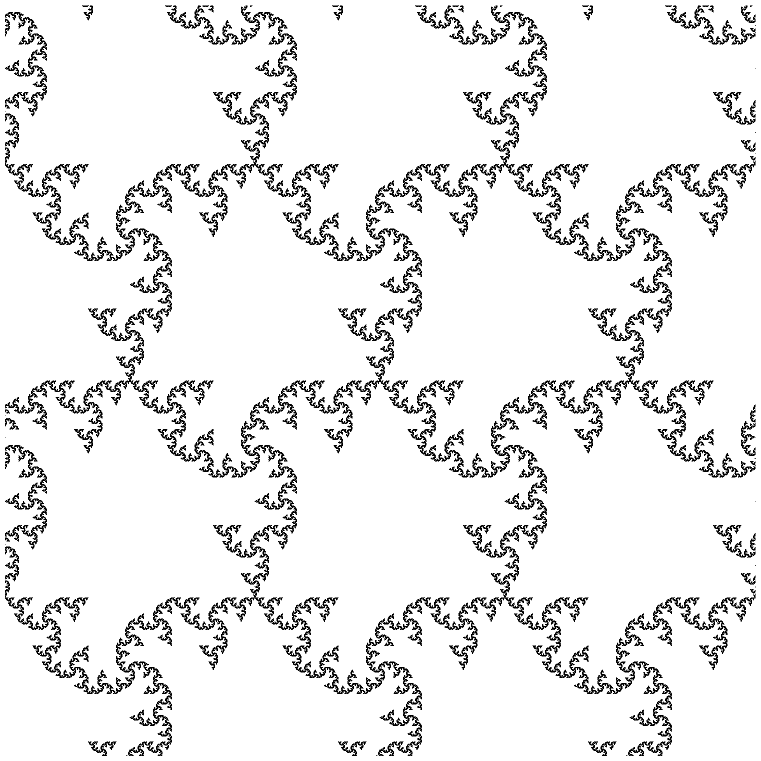}
  \end{center}
  \caption{Covering the plane with triangular tiles}
  \label{fig:plane-tri}
\end{figure}

We can cover the plane with copies of the triangular tile in one size (see figure~\ref{fig:plane-tri}).
This covering is periodic.

\clearpage

\subsubsection{With rhomboidal tiles}

\begin{figure}[!ht]
  \begin{center}
  \includegraphics[height=13cm, keepaspectratio]{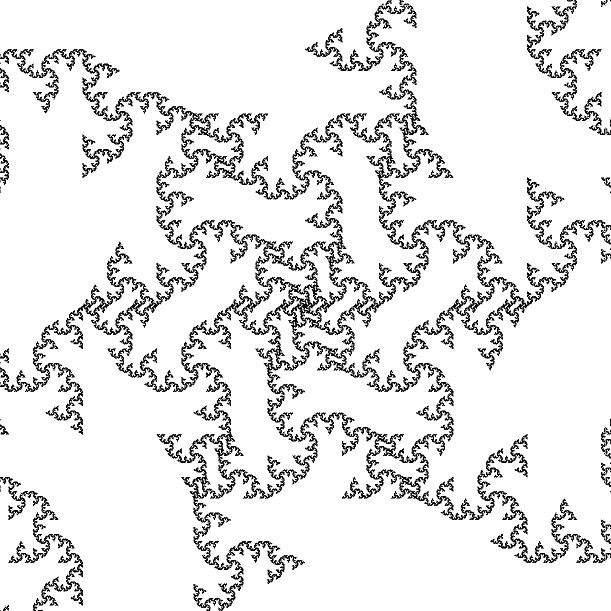}
  \end{center}
  \caption{Covering the plane with rhomboidal tiles}
  \label{fig:plane-rhomb}
\end{figure}

We can cover the plane with copies of the rhomboidal tile scaled
by factors $3^k$ for $k \in \mathbb{Z}$ (see figure~\ref{fig:plane-rhomb}).
This covering has scale invariance.

\clearpage

\subsubsection{With dart tiles}

\begin{figure}[!ht]
  \begin{center}
  \includegraphics[height=13cm, keepaspectratio]{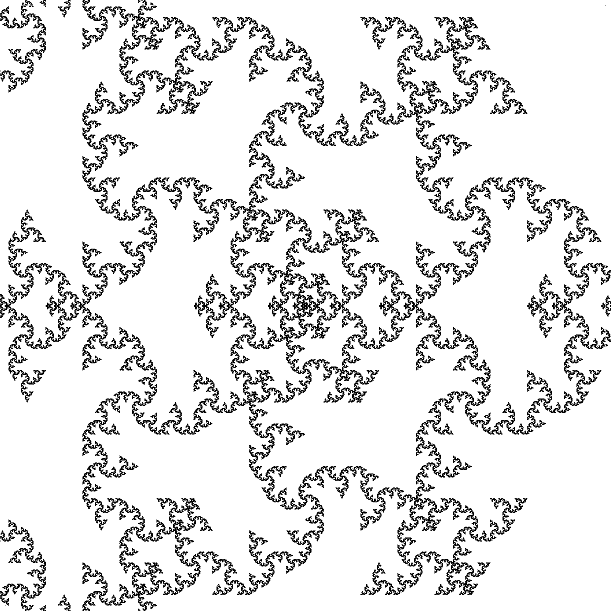}
  \end{center}
  \caption{Covering the plane with dart tiles}
  \label{fig:plane-dart}
\end{figure}

We can cover the plane with copies of the dart tile scaled
by factors $3^k$ for $k \in \mathbb{Z}$ (see figure~\ref{fig:plane-dart}).
This covering has scale invariance.

\section{Acknowledgments}
I would like to express my very great appreciation to J\"org Arndt, Kevin Ryde and Neil Sloane
for their valuable and constructive suggestions.

\end{document}